\documentclass[english]{article} 
\usepackage{mathptmx}
\usepackage{graphicx}
\usepackage{epstopdf}
\usepackage{lineno}
\usepackage{natbib}
\usepackage{enumitem}
\usepackage{amsmath, amssymb, bm}
\usepackage[T1]{fontenc} 
\usepackage[latin9]{inputenc} 
\usepackage[letterpaper]{geometry} 
\usepackage{tikz}
\usepackage{esint}
\usepackage{authblk}
\usepackage{setspace}

\begin{document}

	\title{Characterizing the impact of model error in hydrogeologic time series recovery inverse problems\footnote{Los Alamos National Laboratory technical report: \textbf{LA-UR-16-28825}}}
	
	\author[1]{Scott K. Hansen\footnote{Corresponding author. skh3@lanl.gov}}
	\author[1,2]{Jiachuan He}
	\author[1]{Velimir V. Vesselinov}
	\affil[1]{Computational Earth Science Group (EES-16), Los Alamos National Laboratory}
	\affil[2]{Computational Hydraulics Group (CHG), Institute for Computational Engineering and Sciences (ICES), The University of Texas at Austin}

	\date{February 23, 2017}
	
	\maketitle
	\clearpage

\maketitle
\clearpage

\begin{abstract}
		Hydrogeologic models are commonly over-smoothed relative to reality, owing to the difficulty of obtaining accurate high-resolution information about the subsurface. When used in an inversion context, such models may introduce systematic biases which cannot be encapsulated by an unbiased ``observation noise'' term of the type assumed by standard regularization theory and typical Bayesian formulations. Despite its importance, model error is difficult to encapsulate systematically and is often neglected. Here, model error is considered for a hydrogeologically important class of inverse problems that includes interpretation of hydraulic transients and contaminant source history inference: reconstruction of a time series that has been convolved against a transfer function (i.e., impulse response) that is only approximately known. Using established harmonic theory along with two results established here regarding triangular Toeplitz matrices, upper and lower error bounds are derived for the effect of systematic model error on time series recovery for both well-determined and over-determined inverse problems. A Monte Carlo study of a realistic hydraulic reconstruction problem is presented, and the lower error bound is seen informative about expected behavior. A possible diagnostic criterion for blind transfer function characterization is also uncovered.
\end{abstract}

\section{Introduction}
The effect of systematic model error is being increasingly studied in the hydrological literature, specifically in the context of forward modeling \citep[][]{Gupta2008,Vrugt2008,Lin2012a,Gong2013,Vrugt2013,White2014} and data assimilation \citep[][]{DelGiudice2015a}, with systematic treatments presented by \cite{Refsgaard2006} and \cite{Gupta2012}. Some attention has also been directed explicitly at the effect of model error in inverse modeling \citep{Hansen2016}.

Inverse analyses now form a major part of hydrogeologic practice, being of relevance to all forms of model parametrization, including hydraulic tomography and contaminant source identification. In the inverse problems literature, it is common to assume a perfect model, with all divergence between model prediction and observed data vector attributable to ``noise'' drawn from a symmetric, zero-mean probability distribution function. This theoretical approach underlies classical regularization methods such as Tikhonov and TSVD techniques \citep{Hansen1992}, and is also used for specifying the likelihood function in the Bayesian inversion paradigm \citep{BuiThanh2012}. That the approach of encoding all errors as unbiased parametric uncertainties may not be appropriate in hydrogeologic inverse modeling has been recognized. However in the absence of a paradigm that captures model error in a systematic fashion, the perfect model assumption remains common in practice \citep{Lin2012a,DelGiudice2015}. Thus, it is timely to consider formal analyses of the systematic model errors on inverse-model estimates.

Current approaches to quantifying the effect of model error are typically probabilistic, treating the impact of the model uncertainty on output with Bayesian \citep[][]{Krzysztofowicz1999} or information theoretic \citep[][]{Gong2013} formalisms. The uncertainty about model structure is modeled by parameterizing the model itself as a probability distribution function (pdf) linking inputs and outputs, or as a deterministic numerical model with pdfs defined on its state variables \citep[see summary in][]{Renard2010}. To quantify uncertainty using either approach, some type of Monte Carlo computational analysis is indicated. Since computational exploration cannot proceed in infinite dimensions, it is of course unavoidable that explorations will impose some sort of coarse-grained parametrization which is at best approximately valid.

This paper follows a somewhat different path. The focus is on a specific class of inverse problems commonly faced by hydrogeologists: time series recovery problems with a temporal convolution structure. More concretely, this means the recovery of an input signal from one or more remote output signal measurements, where each output signal is generated by temporal convolution of the (shared) input signal with a (unique) transfer function, which is only approximately known. Hydrogeologic examples of such problems include the inference of hydraulic head history at some location of interest from available time series obtained at remote monitoring wells, and the inference of contaminant source histories from breakthrough curves. 

For these problems, it is shown in Sect. 2 how it is possible to formally decompose the transfer function(s) as well as the input and output signals into generalized Fourier series. Some apparently new results concerning triangular Toeplitz matrices are established. Techniques of matrix analysis are then employed to derive concrete error bounds on the L2 signal reconstruction as a function of the error in dominant components of the transfer function(s). In Sect. 3, a Monte Carlo study of hydraulic inversion is presented which contextualizes the theoretical results shown in Sect. 2 and some empirical observations that go beyond the theoretical work are noted. Section 4 summarizes what has been learned and suggests interesting future research directions.

\section{Derivation of error bounds}

\subsection{Laguerre expansion method}
Consider a linear system in which it is possible to express a transient output signal, $h(\bm{x},t)$, at a location $\bm{x}$, resulting from a transient input signal at location $\bm{0}$, $h(\bm{0},t)$, by means of the convolution $h(\bm{x},t) = b(\bm{x},t)*h(\bm{0},t)$, where $b$ is a Green's function (i.e., transfer function) representing the response to an instantaneous Dirac input signal at the origin (i.e. $h(\bm{0},t)=\delta(t)$). In general, for some fixed $\bm{x}$, the input signal, transfer function, and output signal can be expanded as a generalized Fourier series in a basis of Laguerre functions, ${\phi_n(\cdot)}$ ($n \ge 0$). Each Laguerre function is defined according to the formula:
\begin{equation}
\phi_n(t) = \frac{e^\frac{t}{2}}{n!}\frac{d^n}{dt^n}(e^{-t}t^n),
\end{equation}
and together they form an orthonormal basis on $[0,\infty)$ \citep{Abate1996}. The series expansions are written as follows: 
\begin{eqnarray}
h(\bm{0},t) &=& \sum_n a_n\phi_n\left(\frac{t}{T}\right), \\
b(\bm{x},t) &=& \sum_n b_n\phi_n\left(\frac{t}{T}\right), \\ 
h(\bm{x},t) &=& \sum_n c_n\phi_n\left(\frac{t}{T}\right),
\end{eqnarray}
where $T$ is a characteristic time of the problem, chosen to accelerate convergence.
Let $\mathbf{c}$ be a vector of $N$ Laguerre coefficients, such that its $n$-th entry, $\mathbf{c}_n = c_n$. Similarly, define $\mathbf{a}$ to be a vector of $N$ Laguerre coefficients, such that $\mathbf{a}_n = a_n$. It has been shown \citep{Hansen2009} that, in general, these vectors of coefficients can be related by the matrix operation
\begin{equation}
\mathbf{c=Ba},
\label{eq: exact MIP}
\end{equation}
where $\mathbf{B}$ is the following lower triangular Toeplitz (LTT) matrix:
\begin{equation}
\mathbf{B}= T
\begin{bmatrix} 
b_0			&	0		&	0		& \dots	& 0 \\ 
b_1 - b_0 	& b_0		&	0 		& \dots	& 0 \\ 
b_2 - b_1 	& b_1 - b_0	&	b_0		& \dots	& 0 \\ 
\vdots 		& 	\vdots	& 	\vdots	& \ddots& \vdots  \\ 
b_{N-1} - b_{N-2} 	& b_{N-2} - b_{N-3}		&	b_{N-3} - b_{N-4} 		& \dots	& b_0 \\ 
\end{bmatrix}.
\label{eq: matrix}
\end{equation}
Because $\mathbf{B}$ is a full-rank square matrix, if it is known perfectly then inversion is well defined (although not necessarily numerically stable): 
\begin{equation}
\mathbf{B^{-1}c=a}.
\label{eq: exact inverse MIP}
\end{equation}

\subsection{Results on triangular Toeplitz matrix manipulation}
To continue the analysis of the last section, we need to establish some properties of LTT matrices. First, it is proven that the inverse of a (lower or upper) triangular Toeplitz matrix is itself a (lower or upper) triangular Toeplitz matrix. Second, it is proven that the product of two (both lower or both upper) triangular Toeplitz matrices is similarly a (lower or upper) triangular Toeplitz matrix. Without loss of generality, it is assumed the matrices are LTT in both proofs. 
\subsubsection{Lemma 1: Inversion of triangular Toeplitz matrices}
	Let $\mathbf{M_N}$ be an $N\times N$ LTT matrix, $N$ arbitrary, and let $\mathbf{M_{N}^{-1}}$ be its inverse. It may be shown by induction that $\mathbf{M_{N}^{-1}}$ is LTT. This argument makes repeated use of the following identity for block-triangular matrices \citep[][p. 71]{Bernstein2005}:
	\begin{equation}
	\begin{bmatrix}
	A & 0 \\ C & D
	\end{bmatrix}^{-1} =
	\begin{bmatrix}
	A^{-1} & 0 \\ -D^{-1}CA^{-1} & D^{-1}
	\end{bmatrix},
	\label{eq: block triangle identity}
	\end{equation}
	where $A$, $0$, $C$, and $D$ are compatibly-shaped sub-matrices. 
	
	The base case is trivial: note that that for any $2 \times 2$ LTT matrix, $\mathbf{M_2}$, Eq. (\ref{eq: block triangle identity}) implies directly that $\mathbf{M_2^{-1}}$ is LTT.
	
	For the inductive step, assume that it has been established for $(N-1)\times(N-1)$ LTT matrices that their inverses are LTT. Define $\mathbf{M_{N-1}}$ to be the sub-matrix consisting of the first $N-1$ rows and first $N-1$ columns of an arbitrary LTT matrix, $\mathbf{M_N}$. Note that $\mathbf{M_{N-1}}$ is also LTT, and by our inductive assumption so is $\mathbf{M_{N-1}^{-1}}$. It is valid to apply Eq. (\ref{eq: block triangle identity}) in two different ways. First, make the assignment $A\equiv\mathbf{M_{N-1}}$ and apply Eq. (\ref{eq: block triangle identity}). This  implies that $A^{-1}$ is LTT, and also that $\mathbf{M_{N}^{-1}}$ is lower triangular. This analysis has accounted for all but the $N$-th row of $\mathbf{M_{N}^{-1}}$. To see that the constant descending diagonals continue into the bottom row, note that the sub-matrix consisting of the last $N-1$ rows and last $N-1$ columns of $\mathbf{M_N}$ is also $\mathbf{M_{N-1}}$. Make the assignment $D\equiv\mathbf{M_{N-1}}$ and apply Eq. (\ref{eq: block triangle identity}) again, implying that $D^{-1}$ is LTT, and also that $A^{-1}=D^{-1}$ (Note that these are both $(N-1)\times(N-1)$ matrices which are largely overlapping, and do not participate in the same block partitioning of $\mathbf{M_{N}^{-1}}$). It is thus shown that all descending diagonals of $\mathbf{M_{N}^{-1}}$ are constant (the single element $(\mathbf{M_N^{-1}})_{N1}$ can have any value without affecting this). It has thus been shown that, subject to our inductive assumption, $\mathbf{M_N^{-1}}$ is LTT, for arbitrary LTT $\mathbf{M_N}$. 
	
	By combination of base case and inductive step it follows that if $\mathbf{M_N}$ is an LTT $N \times N$ matrix then so is $\mathbf{M_N^{-1}}$, $\forall\ N \ge 2$. $\square$

\subsubsection{Lemma 2: Multiplication of triangular Toeplitz matrices}
	For any two $N\times N$ matrices $\mathbf{F}$ and $\mathbf{G}$, it is true that the element $(\mathbf{FG})_{ij} = \sum_{k=1}^N\mathbf{F}_{ik}\mathbf{G}_{kj}$. If the matrices are also LTT, it follows that
	\begin{eqnarray}
	\mathbf{F}_{ik} &=& \left\{			
	\begin{matrix}
	0 & i<k \\ 
	\mathbf{f}_{i-k} & i \le k
	\end{matrix}\right. \\
	\mathbf{G}_{kj} &=& \left\{			
	\begin{matrix}
	0 & k<j \\
	\mathbf{g}_{k-j} & k \ge j
	\end{matrix}\right. ,
	\end{eqnarray}
	where $\mathbf{f}_{n}$ and $\mathbf{g}_{n}$ are the elements on the $n$-th diagonal of $\mathbf{F}$ and $\mathbf{G}$, respectively (with the main diagonal having index 0, the sub-diagonal having index 1, and so on). Then it follows that 
	\begin{eqnarray}
	(\mathbf{FG})_{ij} &=& \sum_{k=j}^i\mathbf{f}_{i-k}\mathbf{g}_{k-j}\\
	&=& \sum_{k=0}^{i-j}\mathbf{f}_{(i-j)-k}\mathbf{g}_{k}.
	\label{eq: LTT composition}
	\end{eqnarray} 
	$(\mathbf{FG})_{ij}$ is thus a function only of $i-j$ and is zero for $i<j$. Thus, $\mathbf{FG}$ is LTT. $\square$

\subsection{Effect of imperfect model: single observation location}
Assume perfect knowledge of $h(\bm{x},t)$, but imperfect knowledge of $b(\bm{x},t)$, and a need to infer $h(\bm{0},t)$. The imperfect knowledge of $b$ will lead to an approximate solution $\tilde{h}(\bm{0},t)$, whose Laguerre coefficients lie in vector $\mathbf{\tilde{a}}$. 

In matrix form, this can be written by distinguishing the (unknown) true matrix, $\mathbf{B}$, from the approximate matrix, $\mathbf{\tilde{B}}$, resulting from our imperfect knowledge of the Green's function, $b$. The matrix inverse problem that is being solved is thus
\begin{equation}
\mathbf{c=\tilde{B}\tilde{a}}.
\label{eq: approx MIP}
\end{equation} 
It is also true (by Parseval's theorem) that the squared error of our source history estimate can be expressed in vector form via
\begin{equation}
\int_0^\infty (h(\bm{0},t)-\tilde{h}(\bm{0},t))^2 dt \approx \left \lVert \mathbf{a-\tilde{a}} \right \rVert_2^2,
\end{equation}
with equality in the limit $N\rightarrow\infty$. Although the error introduced by \textit{spectral leakage} (i.e., series truncation) has been recognized as important in some geophysical inversion \citep{Sneider1999}, many transfer functions and input signals are smooth in hydrology, and the approximate equality will be taken to be exact in subsequent analysis. The error analysis can thus be performed in the matrix domain.

Because $\mathbf{\tilde{B}}$ and $\mathbf{B}$ are invertible, there is a unique solution to Eq. (\ref{eq: approx MIP}) and thus:
\begin{equation}
\left \lVert \mathbf{a-\tilde{a}} \right \rVert_2 = \left \lVert \mathbf{(I-\tilde{B}^{-1}B)a} \right \rVert_2.
\label{eq: single error}
\end{equation}
By application of Lemma 1 and Lemma 2, it follows that $\mathbf{I-\tilde{B}^{-1}B}$ is an LTT matrix.
Let $\mathbf{\tilde{b}^{-1}}$, $\mathbf{b}$, and $\mathbf{e}$ be the vectors of coefficients on the diagonals of $\mathbf{\tilde{B}^{-1}}$, $\mathbf{B}$, and  $\mathbf{I-\tilde{B}^{-1}B}$, respectively, indexed as in the proof of Lemma 2. By applying Eq. (\ref{eq: block triangle identity}) with $A$ defined as the upper right $2 \times 2$ sub-matrix of $\mathbf{\tilde{B}}$, it follows immediately that $\mathbf{\tilde{b}^{-1}}_0 = \frac{1}{\tilde{b}_0}$ and $\mathbf{\tilde{b}^{-1}}_1 = - \frac{\tilde{b}_1-\tilde{b}_0}{\tilde{b}_0^2}$. Also, from inspection of Eq. (\ref{eq: matrix}), $\mathbf{b}_0 = b_0$ and $\mathbf{b}_1 = b_1 - b_0$. Then, by applying Eq. (\ref{eq: block triangle identity}) and Eq. (\ref{eq: LTT composition}), it follows that
\begin{eqnarray}
\mathbf{e}_0 &=& 1- \frac{b_0}{\tilde{b}_0}\\
\mathbf{e}_1 &=& \frac{1}{\tilde{b}_0}\left[b_1 -\tilde{b}_1 \left( \frac{b_0}{\tilde{b}_0}\right)\right].
\end{eqnarray}
Note that these elements (like all of $\mathbf{e}$) are zero when $\mathbf{\tilde{B}=B}$.

The following lower error bound follows from consideration of the first element of $\mathbf{(I-\tilde{B}^{-1}B)a}$:
\begin{equation}
|a_0| \left| 1-\frac{b_0}{\tilde{b}_0}\right| \leq \left \lVert \mathbf{a-\tilde{a}} \right \rVert_2,
\end{equation}
In the useful special case in which the input signal is an arbitrary decaying exponential (with rate constant by $1/2T$, noting that $T$ is a free parameter), the only nonzero term of its Laguerre series is $a_0$ and a lower bound on the \textit{relative} error follows immediately:
\begin{equation}
\left| 1-\frac{b_0}{\tilde{b}_0}\right| \leq \frac{\left\lVert \mathbf{a-\tilde{a}} \right\rVert_2}{\left\lVert \mathbf{a} \right\rVert_2}.
\label{eq: lower bound}
\end{equation}
The coefficient $b_0$ is computed
\begin{equation}
b_0 = \int_0^\infty e^{-\frac{t}{2}} b(0,t) dt,
\end{equation}
and similarly for $\tilde{b}_0$. It should be clear that if $b(x,t)$ and $\tilde{b}(x,t)$ have different shapes, particularly if $\tilde{b}(x,t)$ represents transmission through a homogeneous medium, and the true Green's function, $b(x,t)$, is characterizes a medium that is heterogeneous or homogeneous with substantially different properties, then it is possible to have $b_0/\tilde{b}_0 \gg 1$. In such cases, the error due to fitting of the inaccurate model overwhelms the signal and the errors of signal measurement (detection).

It is generally possible to use the approach developed here to generate lower error bounds relating the first $k$ terms of the sequences $\{a_n\}$, $\{b_n\}$, and $\{\tilde{b}_n\}$, for arbitrary $k$, depending on the amount of information available. For instance, for $k=2$: 
\begin{equation}
\left| a_0 \left(1-\frac{b_0}{\tilde{b}_0}\right)\right|^2 + \left| a_0\left(\frac{1}{\tilde{b}_0}\left[b_1 -\tilde{b}_1 \left( \frac{b_0}{\tilde{b}_0}\right)\right] \right) +  a_1 \left(1-\frac{b_0}{\tilde{b}_0}\right) \right|^2 \leq \left \lVert \mathbf{a-\tilde{a}} \right \rVert_2^2.
\label{eq: lower bound2}
\end{equation}
It is also possible to derive an upper bound, which does not depend on $\{a_n\}$, but which requires $k=N$ terms of the other sequences. This is seen in the next section. 

\subsection{Effect of imperfect model: multiple observation locations}
In the case of $M$ monitoring locations, the problem is generally over-determined, and instead of directly computing the inverse, one may define the optimal solution, $\mathbf{\tilde{a}}$, as the one which minimizes the sum of squared residuals at each location, i.e., satisfies the following condition:
\begin{equation}
\sum_{l=1}^M\lVert \mathbf{c_l-\tilde{c_l}} \rVert_2^2= \min_\mathbf{\tilde{a}} \sum_{l=1}^M\lVert \mathbf{c_l-\tilde{c_l}} \rVert_2^2.
\label{eq:srr}
\end{equation}

This problem may be placed in a matrix form by defining the block diagonal matrices
\begin{eqnarray}
\mathbf{B_\otimes} &=&			
\begin{bmatrix} 
\mathbf{B_1} &	0			& \dots		& 0 \\ 
0 			& \mathbf{B_2}	& \dots		& 0 \\ 
\vdots 		& 	\vdots		& \ddots	& \vdots  \\ 
0			& 0				& \dots		& \mathbf{B_M}
\end{bmatrix},\\
\mathbf{\tilde{B}_\otimes} &=&
\begin{bmatrix} 
\mathbf{\tilde{B}_1} &	0			& \dots		& 0 \\ 
0 			& \mathbf{\tilde{B}_2}	& \dots		& 0 \\ 
\vdots 		& 	\vdots		& \ddots	& \vdots  \\ 
0			& 0				& \dots		& \mathbf{\tilde{B}_M}
\end{bmatrix}.
\end{eqnarray}
It is also useful to define the following block-columnar matrix of $M$, $N\times N$ identity matrices:
\begin{equation}
\mathbf{D} =
\begin{bmatrix} 
\mathbf{I_N}\\ 
\mathbf{I_N} \\ 
\vdots 		 \\ 
\mathbf{I_N}
\end{bmatrix}.			
\end{equation}
Using this notation, Eq. (\ref{eq:srr}) can be represented as:\textsl{}
\begin{equation}
\sum_{l=1}^M\lVert \mathbf{c_l-\tilde{c_l}} \rVert_2^2 = \mathbf{a^T D^T B_\otimes^T B_\otimes D a} -2\mathbf{a^T D^T B_\otimes^T \tilde{B}_\otimes D \tilde{a} + \tilde{a}^T D^T \tilde{B}_\otimes^T \tilde{B}_\otimes D \tilde{a}}.
\end{equation}
Differentiating,
\begin{eqnarray}
\frac{d}{d\mathbf{\tilde{a}}}\sum_{l=1}^M\lVert \mathbf{c_l-\tilde{c_l}} \rVert_2^2
&=&\frac{d}{d\mathbf{\tilde{a}}} [-2\mathbf{a^T D^T B_\otimes^T \tilde{B}_\otimes D \tilde{a} + \tilde{a}^T D^T \tilde{B}_\otimes^T \tilde{B}_\otimes D \tilde{a}}]\\
&=&-2\mathbf{a^T D^T B_\otimes^T \tilde{B}_\otimes D} + 2 \mathbf{\tilde{a}^T D^T \tilde{B}_\otimes^T \tilde{B}_\otimes D}
\end{eqnarray}
The optimal solution will be when this quantity equals zero, which is satisfied when
\begin{equation}
\mathbf{ D \tilde{a} = \tilde{B}_\otimes^{-1} B_\otimes D a },
\end{equation}
or
\begin{equation}
\mathbf{ \tilde{a}} = \frac{1}{M} \mathbf{D^T\tilde{B}_\otimes^{-1} B_\otimes D a }.
\end{equation}
It is clear that
\begin{eqnarray}
\left\lVert \mathbf{a-\tilde{a}} \right\rVert_2 &=& \frac{1}{M}\left\lVert \mathbf{Da-D\tilde{a}} \right\rVert_2\\
&=& \frac{1}{M}\left\lVert \mathbf{(I_{MN}-\tilde{B}_\otimes^{-1}B_\otimes) Da} \right\rVert_2. \label{eq: kron error}
\end{eqnarray}
By inspection of the last equality, it is apparent that
\begin{eqnarray}
\left\lVert \mathbf{a-\tilde{a}} \right\rVert_2 =\frac{1}{M}  \sum_{l=1}^M \left\lVert \mathbf{(I_{N}-\tilde{B}_l^{-1}B_l) a} \right\rVert_2. 
\label{eq: series error}
\end{eqnarray}
This is simply the average of the individual model errors if only a single monitoring location were to be used (see Eq. (\ref{eq: single error})), implying that the error bound theory developed above can be carried over straightforwardly. If $\mathbf{a}$ is viewed as a random variable, the expected error is also not reduced by incorporating additional measurements, unlike the scenario of uncorrelated random noise.

One can also derive an upper bound on the relative error of $\mathbf{\tilde{a}}$ from Eq. (\ref{eq: kron error}) by noting that the H\"{o}lder norm $\lVert\cdot\rVert_2$ is submultiplicative, and that $\frac{1}{M}\lVert\mathbf{Da}\rVert_2 = \lVert\mathbf{a}\rVert_2$. From this, it follows immediately that
\begin{equation}
\frac{\left\lVert \mathbf{a-\tilde{a}} \right\rVert_2}{\left\lVert \mathbf{a} \right\rVert_2} \leq \left\lVert \mathbf{I_{MN}-\tilde{B}_\otimes^{-1}B_\otimes}\right\rVert_2.
\end{equation}
Naturally, this applies to single measurement location reconstruction as a special case. Practically, this depends on the full sets of coefficients $\{b_n\}$ and $\{\tilde{b}_n\}$, which is a greater information demand than for the lower bounds, which only involved relationships of the dominant components.

\section{Monte Carlo study: reconstruction of hydraulic transients}

In this section, the inference of a hydraulic head transient history along an aquifer boundary (which might be interpreted as a river stage transient, where the river is in a hydraulic connection with the aquifer) based on a time series of measurements made at a single nearby groundwater monitoring well is considered. This represents both an application of the above ideas, and also a study of independent interest. 

\subsection{Procedure}

It is assumed here that the specific storage is known and spatially uniform, and that the log hydraulic conductivity is defined by a multi-Gaussian spatially random field whose mean is known, but which is otherwise unknown. Assuming flow is described by the groundwater flow equation on this heterogeneous conductivity field, $b(\bm{x},t)$ is determined as the head history at a fixed location, $\bm{x}$. A natural interpretive model, $\tilde{b}(\bm{x},t)$, is selected: the same groundwater flow equation, but solved on a homogeneous conductivity field that is everywhere equal to the mean of the true log hydraulic conductivity field. 

The study is then performed according to the following basic procedure:
First, a true, exponentially decaying, transient in the river stage on the aquifer domain boundary is specified, along with the location of a monitoring well at which a time series of measurements is to be made. The accuracy of reconstruction of the river stage transient from the transient at the monitoring well is studied, given an overly smooth model of the subsurface. For simplicity, the free parameter, $T$, is selected so that the Laguerre decomposition of the true transient is the vector $\mathbf{a}=<1,0,0,\dots,0>$.
Next, 500 two-dimensional subsurface realizations are generated with different heterogeneous log-hydraulic conductivity fields, all of which have the same multi-Gaussian statistical correlation structure and geometric mean conductivity.
Subsequently, using finite element analysis, head time series are computed at the monitoring well for each of the 500 subsurface realizations resulting from a Dirac head impulse at the river stage. Each impulse response (Green's function) is decomposed as a vector of Laguerre coefficients, $\mathbf{b}$.
Again, using finite element analysis, the impulse response at the well is computed, but assuming a uniform hydraulic conductivity field with the same geometric mean hydraulic conductivity as used in each of the heterogeneous realizations. This impulse response is decomposed as a vector of Laguerre coefficients, $\mathbf{\tilde{b}}$.
Finally, for each realization, the reconstruction error, $\left \lVert \mathbf{(I-\tilde{B}^{-1}B)a} \right \rVert_2$, is computed and compared with the analytical lower bound in Eq. (\ref{eq: lower bound}). Statistics about this quantity are tabulated so that its relationship to qualitative features of the inverse model discrepancy may be studied.

Both the true solution and the interpretive model are described by the following system of equations:
\begin{eqnarray}
S_s \frac{\partial h(\bm{x},t)}{\partial t} + \nabla \cdot \bm{q}(\bm{x},t) &=& 0, \quad \text{$\bm{x} \in \bm{D}$} \label{eq: gwfe}\\
\bm{q}(\bm{x},t) &=& -K(\bm{x}) \nabla h(\bm{x},t), \quad \text{$\bm{x} \in \bm{D}$} \label{eq: darcy}
\end{eqnarray}
solved subject to the initial and boundary conditions
\begin{eqnarray}
h(\bm{x},0) &=& 0, \quad \text{$\bm{x} \in \bm{D}$}\\
h(\bm{x},t) &=& c(t),  \quad \text{$\bm{x} \in \Gamma_L$}\\
h(\bm{x},t) &=& 0,  \quad \text{$\bm{x} \in \Gamma_R$}\\
\bm{q}(\bm{x},t) \cdot \bm{n}(\bm{x}) &=& 0, \quad \text{$\bm{x} \in \Gamma_S$} \label{eq: model}
\end{eqnarray}
where the only difference enters due to different choices for $K(\bm{x})$. In the above equations, $\bm{x} \: [\mathrm{L}]$ represents the location, $t \: [\mathrm{T}]$  represents time, $S_s \: [\mathrm{L^{-1}}]$ represents specific storage, $h \: [\mathrm{L}]$ represents hydraulic head, $\bm{q} \: [\mathrm{LT^{-1}}]$ represents groundwater flux, $K \: [\mathrm{LT^{-1}}]$ represents hydraulic conductivity, and $\bm{n} \: [1]$ is the outward-facing unit normal vector. For vector quantities, the units reported are for each of their components. 

\begin{figure}[h]
	\centering
	\includegraphics{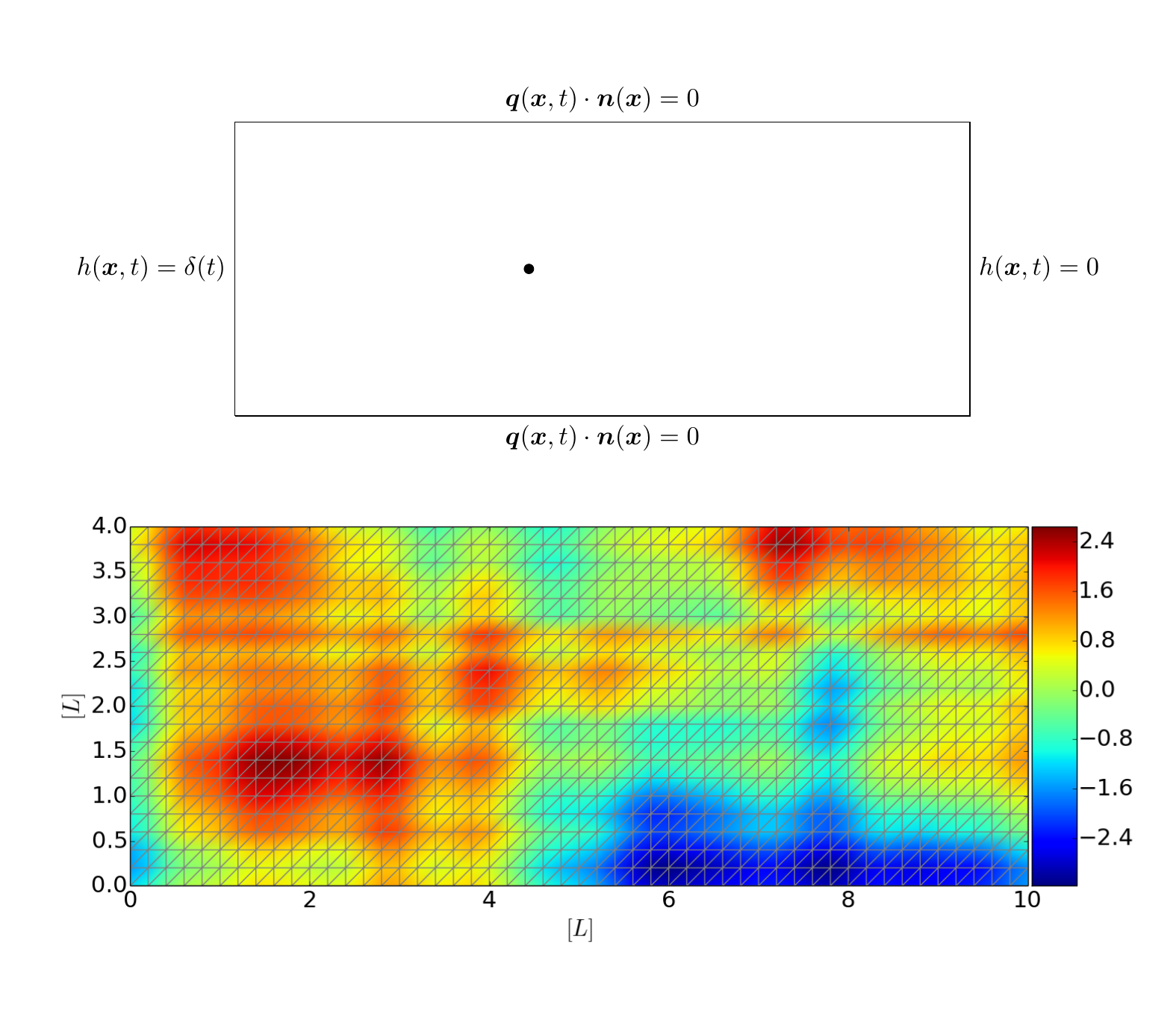}
	\caption{Top: the flow domain and the observation location (black dot). Bottom: example of a realization of $\ln [K(\bm{x})]$ field along with the mesh used to numerically solve Eq. \eqref{eq: model}} 
	\label{flow_domain}
\end{figure}

More concretely, a two-dimensional model of saturated flow in a heterogeneous porous medium is defined over the rectangular domain $\bm{D}=(0,L_1) \times (0,L_2)$, where $L_1=10\ [L]$ and $L_2=4\ [L]$ ($L$ is any consistent length unit), and specific storage, $S_s=1$  (Fig. \ref{flow_domain}). $\Gamma_L$ represents the left boundary of $\bm{D}$ (at $x=0$), $\Gamma_R$ represents the right boundary of $\bm{D}$ (at $x=10$), and $\Gamma_S$ represents the union of the other two sides of $\bm{D}$ (at $y=0$ and $y=4$, respectively). Let $Y(\bm{x}, \omega) = \ln [K(\bm{x}, \omega)]$ be a random field, where $\omega$ belongs to the space of random events $\bm{\Omega}$. Assuming $Y(\bm{x}, \omega)$ is Gaussian with zero mean and a separable exponential covariance function,
\begin{equation}
C(\bm{x_1},\bm{x_2}) = C(x_1, y_1; x_2, y_2) = \sigma^{2}_Y \exp\left[-\frac{\left|x_1 - x_2\right|}{\eta_1} - \frac{\left|y_1 - y_2\right|}{\eta_2}\right],
\end{equation}
where $\sigma^{2}_Y=2$, $\eta_1=4$ and $\eta_2=2$ are the variance and the correlation lengths of the random field. 

\begin{figure}[h]
	\centering
	\includegraphics{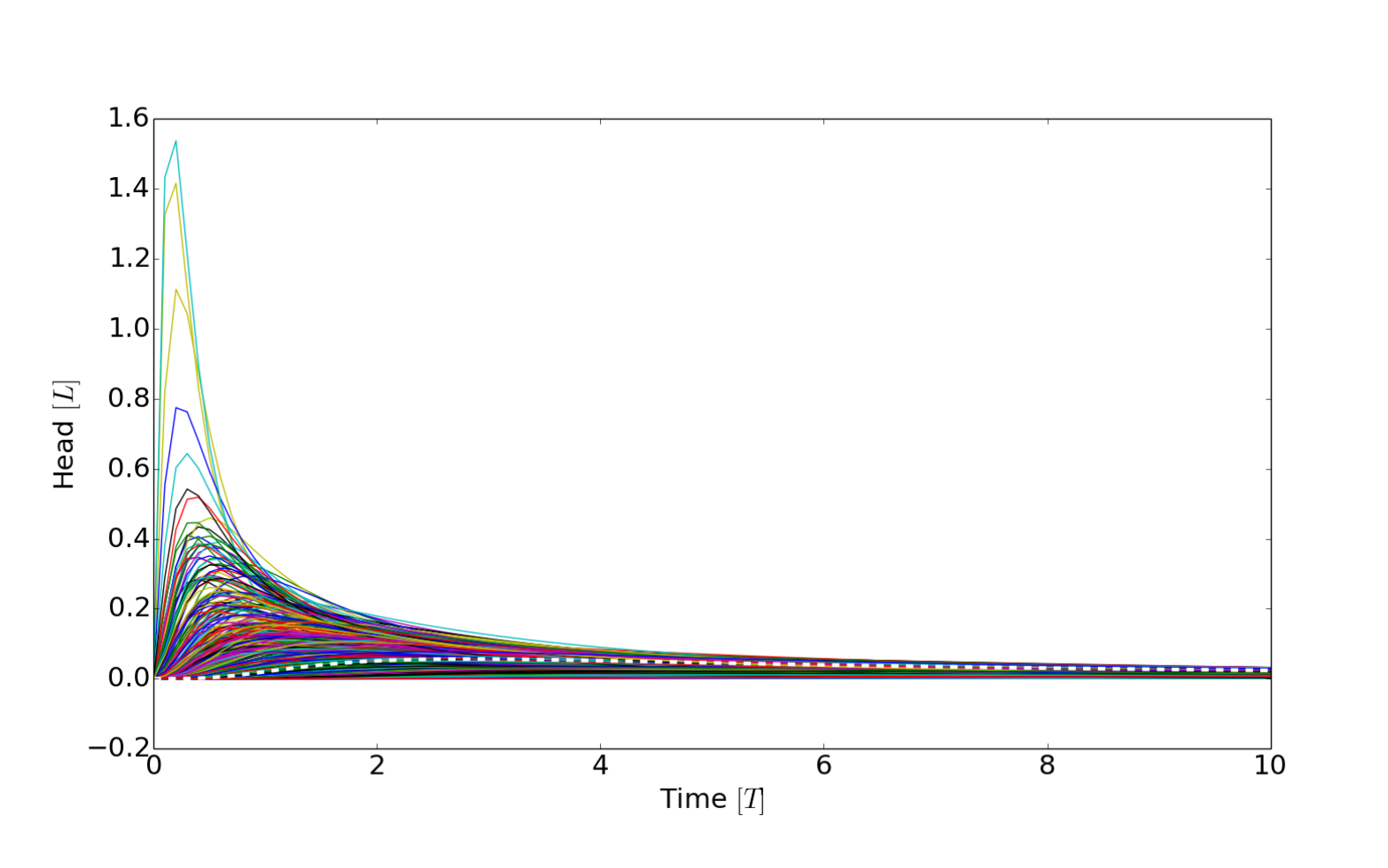}
	\caption{Hydraulic head transient responses at the point $(x,y)=(4,2)$ due to an impulse at the boundary $x=0$. Impulse responses---i.e., $b(4,2,t)$---for 500 random realizations of the hydraulic conductivity field are shown (solid, colored lines), along with the response from the interpretation model---i.e., $\tilde{b}(4,2,t)$---with spatially uniform hydraulic conductivity, $K(\bm{x})=1\  [LT^{-1}]$ (dashed, white line)} 
	\label{observation}
\end{figure}

\begin{figure}[h]
	\centering
	\includegraphics{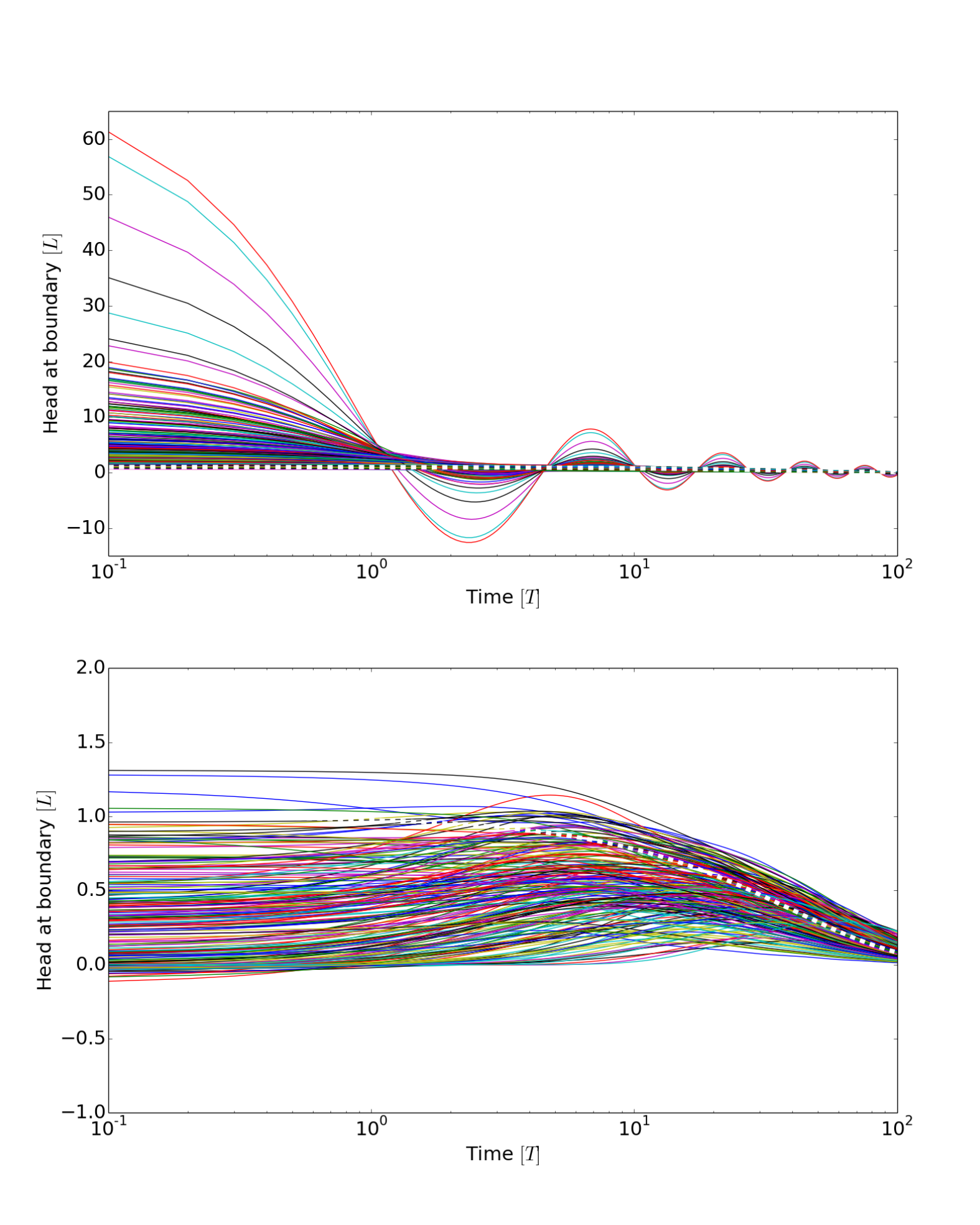}
	\caption{Reconstructions of the boundary condition resulting from use of the interpretation model to 500 simulated impulse responses and the simulated response from the interpretation model (dashed white line). Responses are partitioned according to whether true peak hydraulic head is earlier (top axes) than assumed by the interpretive model, or later (bottom axes)}
	\label{reconstruction}
\end{figure}

For the Monte Carlo study, a set of 500 realizations of $b(t)$ is generated by setting $c(t) = \delta(t)$, generating 500 $\ln [K(\bm{x})]$ fields (Fig. \ref{flow_domain}, bottom), using a 100-term truncated Karhunen-Lo\`eve expansions (KLE) to represent the field as weighted sums of predefined spatially variable orthonormal functions \citep{Zhang2004a}, and solving Eqs. (\ref{eq: gwfe}-\ref{eq: model}) on each.
The numerical solution is evaluated using the FEniCS \citep{Logg2012} package to discretize Eqs. (\ref{eq: gwfe}-\ref{eq: darcy}), using finite element methods in space and an implicit Euler method in time. Simulated hydraulic head data is recorded at the point (4, 2) for each realization of	$K(\bm{x})$ (Fig. \ref{observation}), and the 50-term LEM is used along with $T=100\ [T]$ to reconstruct the time series of hydraulic head values on the left boundary (Fig. \ref{reconstruction}). The same procedure is followed to generate the interpretive model, $\tilde{b}(t)$ except $K(\bm{x}) = 1 \: \forall \bm{x}$ (i.e., $Y(\bm{x}) = 0 \: \forall \bm{x}$) is employed.

\subsection{Discussion of results}

Given that the random hydraulic conductivity fields chosen for the forward modeling were only moderately heterogeneous, and choice of a spatially uniform interpretive model is a natural response to unresolved heterogeneity, the wide array of possible reconstructions is notable. In particular, a bifurcation of the response classes was noted based on whether the peaks of the forward model impulse response, $b$, preceded or lagged that of the interpretive model, $\tilde{b}$. In case in which the interpretive model predicted a faster response than existed in reality, the reconstruction of the decaying exponential boundary condition was typically smooth, with its peak at a time significantly greater than zero: the delayed reconstruction of the boundary condition compensated for the over-rapid model. However, because signal causality is enforced, the model cannot respond to an earlier-than-anticipated arrival with a non-zero signal at negative time. Instead, the optimal reconstruction features a large peak at time zero, followed by decaying corrective oscillations. This bifurcation of behavior is potentially useful as a model diagnostic tool that does not require any a priori knowledge of the true model (other than that it possesses a unimodal structure): multiple candidate interpretive models could be tested with peaks at different locations, and the true peak location pinpointed by the disappearance of the spurious oscillations.

\begin{figure}[h]
	\centering
	\includegraphics{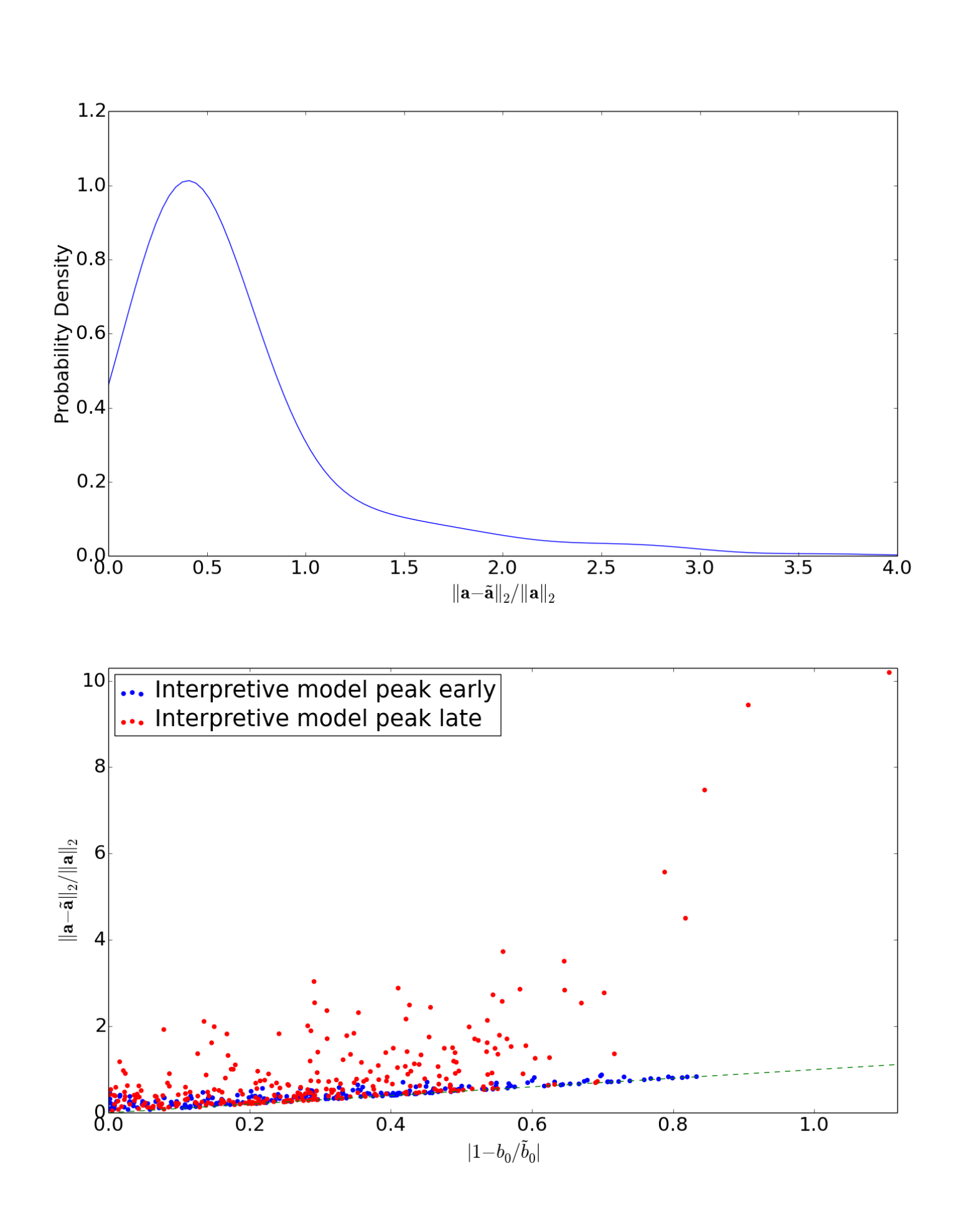}
	\caption{Top: empirical pdf of normalized L2 estimation error of $h(x=0,y,t)$ from the ensemble of 500 realizations. Bottom: scatter plot of relative estimation error against lower error bound}
	\label{fig: error summary}
\end{figure}

In Fig. \ref{fig: error summary} (top), the empirical pdf for the L2 error, normalized by the L2 norm of the signal being reconstructed is shown. It is apparent that even for moderate heterogeneity, reconstruction error on the same magnitude as the signal itself is to be expected. In Fig. \ref{fig: error summary} (bottom), the L2 error of approximation in the reconstruction of the boundary condition for each of the 500 realizations is plotted against the zero-order lower error bound in Eq. \eqref{eq: lower bound}.

\section{Summary and conclusion}
Systematic model error was considered in the context of inverse problems in systems whose output signal is determined by convolution of an input signal with a transfer function, or impulse response, which describes system behavior. Using a generalized Fourier series expansion in Laguerre basis functions, it was possible to translate the signal reconstruction inverse problem into a matrix inverse problem whose structure may be analyzed using some classic, and some apparently new, results in matrix algebra. It is thus seen possible to place upper and lower bounds on the L2 signal reconstruction error as in terms of the transfer function infidelity.

The inverse problem of recovering river level history from remote measurements at a well, which has convolution structure, was chosen for Monte Carlo study. Forward predictions were generated by solving the groundwater flow equation on a mildly heterogeneous hydraulic conductivity field, and these were interpreted using the groundwater flow equation, assuming a homogeneous hydraulic conductivity field. The L2 reconstruction error for the river level was established for all realizations, and this error was compared with the error bound developed above. The simple lower bound derived here (Sect. 2) was found to be informative regarding the reconstruction error in the specific realizations. A qualitative bifurcation in the reconstructed signal was discovered, depending on the location of the peak of the interpretation model transfer function relative to that of the true model. Looking forward, this may prove to be a useful tool for transfer function identification.

The Laguerre expansion approach, because of its high degree of structure, relative simplicity and computational efficiency, may also prove to be a profitable foundation for further analysis of systematic model error. The systematic model errors are commonly ignored in the theoretical and practical inverse-model analyses. The matrix transformation of the inverse problem developed here can be also applied for other problems of interest such as groundwater contaminant transport, propagation of low-frequency seismic waves, heat flow, infectious disease transmission, population dynamics, spreading chemical/biochemical substances in atmosphere, and many others.

\section*{Acknowledgements}
	The authors acknowledge the support of the LANL Environmental Programs. All data is synthetic; the authors maintain an archive of codes and simulation output employed in the paper.

\pagebreak
\pagebreak
\pagebreak
\pagebreak

\end{document}